\documentclass[12pt, reqno]{amsart}
\usepackage{amsmath, amsthm, amscd, amsfonts, amssymb, graphicx, color}
\usepackage[bookmarksnumbered, colorlinks, plainpages]{hyperref}
\hypersetup{colorlinks=true,linkcolor=red, anchorcolor=green, citecolor=cyan, urlcolor=red, filecolor=magenta, pdftoolbar=true}
\usepackage{enumitem}
\newtheorem{theorem}{Theorem}[section]
\newtheorem{lemma}[theorem]{Lemma}

\theoremstyle{definition}

\theoremstyle{remark}
\newtheorem{remark}[theorem]{Remark}
\numberwithin{equation}{section}

\begin{document}

\setcounter{page}{1}

\title[Elliptic gradient estimates and Liouville theorems]{Elliptic gradient estimates and Liouville theorems for a weighted nonlinear parabolic equation}

\author[A. Abolarinwa]{Abimbola Abolarinwa}

\address{Department of Mathematics and Statistics,
 Osun State College of Technology, P. M. B. 1011, Esa-Oke, Osun State,
Nigeria.}
\email{\textcolor[rgb]{0.00,0.00,0.84}{A.Abolarinwa1@gmail.com}}





\date{July 28, 2017}

\begin{abstract}
Let $(M^N, g, e^{-f}dv)$ be a complete smooth metric measure space with $\infty$-Bakry-\'Emery Ricci tensor bounded from below. We derive elliptic  gradient estimates for positive  solutions of a  weighted  nonlinear parabolic equation
\begin{align*}
\displaystyle \Big(\Delta_f - \frac{\partial}{\partial t}\Big) u(x,t) +q(x,t)u^\alpha(x,t)  = 0,
\end{align*}
where $(x,t) \in M^N \times (-\infty, \infty)$ and $\alpha$ is an arbitrary constant.
As Applications we prove a Liouville-type theorem for positive ancient solutions and Harnack-type inequalities for positive bounded solutions.
\end{abstract}
 \maketitle

\section{Introduction}
\subsection{Background}
Recently, there have been many interesting  results relating to parabolic and elliptic gradient estimates, Harnack inequalities  and Liouville-type theorems on either Riemannian manifolds or smooth metric measure spaces. This is due to the fact that these estimates are now fundamental tools in Geometric Analysis and PDEs. Historically, gradient  and Harnack estimates  for parabolic equations on manifold originated in Li and Yau \cite{[LY86]} where they extended the work in \cite{[CY75]}. Then Hamilton \cite{[Ha93]} proved an elliptic type gradient estimate for the heat equation. But this Hamilton-type of estimates is a global result which requires the heat equation defined
on closed manifolds. Souplet and Zhang \cite{[SZ06]} later proved  a localized  version of Hamilton-type gradient estimate by combining  Li-Yau's Harnack
inequality \cite{[LY86]} and Hamilton’s gradient estimate \cite{[Ha93]}. See for examples \cite{[BCP10],[Br13],[CC09],[MH17],[Ru10],[Wa15],[Wu10b],[Wu14],[YY08],[Zhu11],[Zhu13],[Zhu16]} for many more interesting results and applications in various settings. In particular, Brighton \cite{[Br13]} proved an elliptic gradient estimate for positive weighted-harmonic functions by applying Yau's idea to function $u^\epsilon \ (0 < \ \epsilon < 1)$ instead of $\ln u$ used in \cite{[Y75]}, and hence obtained a Liouville theorem for positive bounded weighted harmonic functions with nonnegative $\infty$-Bakry-Emery Ricci tensor.

This paper is an extention of \cite{[Wu16]}. In that paper, the  author proved elliptic gradient estimates and Liouville type theorem for
positive solutions to a nonlinear parabolic equation \begin{align}
\displaystyle \Big(\Delta_f - \frac{\partial}{\partial t}\Big) u(x,t) +a u(x,t) \ln u(x,t)  = 0, \ \ \ a \in \mathbb{R}
\end{align}
on complete smooth metric measure spaces with $m$-Bakry- \'Emery Ricci tensor bounded below. 
In the present paper, we use a similar approach to derive localised elliptic (space only) gradient estimates for positive  solutions to the  weighted  nonlinear parabolic equation
\begin{align}\label{eq11}
\displaystyle \Big(\Delta_f - \frac{\partial}{\partial t}\Big) u(x,t) +q(x,t)u^\alpha(x,t)  = 0,
\end{align}
where $ \alpha \in \mathbb{R}$. The function $q(x,t)$ is a space-time smooth function atleast $C^1$ in space and $C^0$ in time.  If $q(x,t) = 0$, then the nonlinear equation (\ref{eq11}) reduces to a weighted heat equation which was studied in \cite{[Wu15]}. It is well-known that all solutions to Cauchy problem for the weighted heat equation exist for all time. For the applications of our gradient estimates, we prove
parabolic Liouville properties and Harnack inequalities for positive ancient solutions to (\ref{eq11}) under the assumption that $\infty$-Bakry-\'Emery tensor is bounded below.  Notice also that (\ref{eq11}) is a weighted version of
\begin{align}
\displaystyle \Big(\Delta - \frac{\partial}{\partial t}\Big) u(x,t) + q(x,t)u^\alpha(x,t)  = 0, \ \ \ \alpha \geq 1
\end{align}
considered by Zhu in \cite{[Zhu16]}. But in our case $\alpha$ is an arbitrary constant rather than only being greater than $1$.

\subsection{Smooth metric measure spaces} A smooth metric measure space is denoted by the triple $(M^N, g, e^{-f}dv)$,  where $(M^N, g)$ is an $N$-dimensional  complete manifold with the Riemannian metric tensor $g$, volume element  $dv$ and $f$ is a $C^\infty$ real-valued function on $M$. Smooth metric measure spaces are naturally endowed with analogue of Laplace-Beltrami operator, called weighted Laplacian and analogue of Ricci tensor, called Bakry-\'Emery tensor.
The weighted Laplacian defined by
$$\Delta_f := \Delta - \langle \nabla f, \nabla \cdot \rangle,$$
where $\Delta$ is the Laplace-Beltrami operator, is symmetric and self-adjoint with respect to the weighted measure $e^{-f}dv$.
The $m$-Bakry-\'Emery tensor is defined by 
$$Ric_f^m := Ric + \nabla^2 f - \frac{1}{m} df \otimes df$$
for some constant $m>0$, where $Ric$ is the Ricci tensor of the manifold and $\nabla^2$ is the Hessian  with respect to the metric $g$. When $m$ is infinite we have the $\infty$-Bakry-\'Emery tensor
$$ Ric_f =  \lim_{m \to \infty} {Ric_f^m} := Ric + \nabla^2 f.$$
This tensor is related to the gradient Ricci  soliton
$$Ric_f = \lambda g$$
where $\lambda$ is a real constant. A Ricci soliton is said to be shrinking, steady or expanding depending on whether $\lambda$ is positive, zero or negative respectively. Ricci solitons play an imporatnt role in the theory of singularities for the Ricci flow \cite{[Ha95]} (see \cite{[CaH]} for a recent survey on Ricci solitons).
The weighted Laplacian and the Bakry-\'Emery tensor are related by Bochner formula 
\begin{align}\label{eq13}
\frac{1}{2}\Delta_f |\nabla u|^2 = |\nabla^2 u|^2 + \langle \nabla \Delta_fu, \nabla u\rangle + Ric_f(\nabla u, \nabla u).
\end{align}

Since $\Delta_f$ and $Ric_f$ are natural extension of Laplace-Beltrami operator and Ricci tensor respectively, it is not unexpected that many geometric and topological results for Riemannian manifolds could be extended to smooth metric measure spaces, see for instance \cite{[BQ00],[Li05]}.

\subsection{Motivations}
The motivations for this work come from geometric and physical applications of (\ref{eq11}). For instance,  the authors in \cite{[BRS]} show that if $f$ is a constant then the equation
\begin{align}\label{eqa}
\displaystyle \Delta u +q(x)u^\alpha + b(x)u  = 0
\end{align}
which is a non-weighted static version of (\ref{eq11}) (when $b(x) \equiv 0$), is equivalent to Yamabe problem on noncompact Riemannian manifold. Clearly, setting $\widetilde{g}=u^{4/n-2}g, \ u >0$, then for $\mathcal{R}(x)$, the scalar curvature of $g$ and $K(x) \in C^\infty(M)$, the scalar curvature of $\widetilde{g}$, we have the relation
\begin{align}\label{eqb}
\displaystyle \Delta u - \frac{n-2}{4(n-1)}\mathcal{R}(x)u + \frac{n-2}{4(n-1)} K(x) u^{\frac{n+2}{n-2}} = 0,
\end{align}
which is of the form (\ref{eqa}). Yamabe problem demands the existence of a positive everywhere defined solution of (\ref{eqb}). Indeed, the existence and uniqueness of such solution depends on the geometry of the underlying manifold. Thus, $g$ can be pointwise conformally deformed to a complete metric $\widetilde{g}$ of a scalar curvature $K(x)$. For further discussions on existence, uniqueness and a priori estimates of (\ref{eqa}) (resp. Yamabe-type equation) see \cite{[MRS]}.

On the other hand, the static form of (\ref{eq11}) for a special $\alpha$ is related to the Euler-Lagrange equation for the weighted Yamabe quotient on compact smooth metric spaces
\begin{align}\label{eqc}
\displaystyle \Delta_f u - \frac{m+n-2}{4(m+n-1)}\mathcal{R}^m_f u -c_1 u^{\frac{m+n}{m+n-2}}  e^{\frac{f}{m}} + c_2 u^{\frac{m+n+2}{m+n-2}} = 0, \ \ \ m>0
\end{align}
where $\mathcal{R}^m_f$ is the weighted scalar curvature defined by
$$\mathcal{R}^m_f:= \mathcal{R} + 2 \Delta f - \frac{m+1}{m}|\nabla f|^2.$$
Thus, the weighted volume can be conformally deformed as in \cite{[Cas2]}. In fact, setting
$$(M, \widetilde{g}, e^{-\widetilde{f}} d\widetilde{v},m) = (M, e^{\frac{2\sigma}{m+n-2}}g, e^{\frac{(m+n)\sigma}{m+n-2}}e^{-f}dv)$$
for some $\sigma \in C^\infty(M)$, then the weighted Yamabe quotient is conformally invariant, i.e., $\widetilde{Q}(u) = Q(e^{\sigma/2}u)$, where $Q : C^\infty(M) \to \mathbb{R}$ is defined by the functional
$$Q(u):= \frac{\Big(\int|\nabla u|^2 + \frac{m+n-2}{4(m+n-1)}\mathcal{R}^m_f u^2 \Big)\Big( \int |u|^{\frac{2(m+n-1)}{m+n-2}} e^{\frac{f}{m}}\Big)^{\frac{2m}{n}}}{\Big(\int |u|^{\frac{2(m+n)}{m+n-2}}\Big)^{\frac{2m+n-2}{n}},}$$
where all integrals are with respect to the weighted measure. The infimum of the above functional for all $u \in W^{1,2}(M,g, e^{-f}dv)$ is called the weighted Yamabe constant.
 Case also shows in \cite{[Cas1]} that Yamabe-type problem on $(M,g, e^{-f}dv)$ interpolates between Yamabe problem and the problem of finding minimizers for Perelman's $\nu$-entropy. Meanwhile, it is well known that Yamabe constant and Perelman's $\nu$-entropy are remarkably related to  Sobolev and logarithmic Sobolev inequalities. Interested readers can check for \cite{[DD],[Le]} and \cite{[Cas1]}. We do hope that the gradient estimate for equation (\ref{eq11} may be useful for tackling the Yamabe problem of smooth metric measure spaces.

Suppose $f$ is a constant function on $M$, a physical application of the term $q(x,t)u^\alpha$  of (\ref{eq11}) may be seen if $u=u(x,t)$ is considered as a population density. Here the nonlinearity in $u^\alpha$ could be interpreted as intraspecies interraction like competition or inhibition, while its product with $q(x,t)$ is a spatial relation in the form of interraction with enviroment. If such interraction is time-independent, then $q(x,t)$ would be replaced by $q(x)$.

\subsection{Main results}
This paper aims majorly at obtaining elliptic type gradient estimates for positive solutions of (\ref{eq11}).
Precisely, Let 
$$\mathcal{Q}_{R,T} \equiv B(x_0,R) \times [t_0-T, t_0] \subset M \times (-\infty, \infty),$$
 where $B(x_0,R)$ is a ball of radius $R>0$ centred at $x_0$, and $t_0 \in \mathbb{R}$, $T>0$, we have

\begin{theorem}\label{thm11}
Let $(M^N, g, e^{-f}dv)$ be an $N$-dimensional  complete smooth metric measure space with $Ric_f \geq -(N-1)K$ for some $K \geq 0$. Fix $x_0 \in M$ and $R \geq 2$. 
Suppose that $u(x,t)$ is a positive solution to (\ref{eq11}) in $\mathcal{Q}_{R,T}$, $T>0$. Suppose further that $u(x,t) \leq D$ for some constant $D$ in $\mathcal{Q}_{R,T}$ and $\beta := \sup_{(x,t) \in \mathcal{Q}_{R/2,T}}|h| +1 $, where $h = \ln u/D$. Then, for all $(x,t) \in \mathcal{Q}_{R/2,T}$ with $t \neq t_0-T$, there exists a constant $C(\delta)$ depending on $N$ and $\delta$  such that 
\begin{enumerate}
\item if $\alpha \geq 1$
\begin{equation}\label{eq12}
\left. \begin{array}{l}
\displaystyle \frac{|\nabla u(x,t)|}{u(x,t)} \leq  C(\delta)  \Big( \sqrt{\frac{1+|\mu|}{R}} + \frac{1}{\sqrt{t-(t_0-T))}} +  \sqrt{K}   +  \sqrt{\alpha} D^{\frac{1}{2}(\alpha-1)}\|q^+\|^{1/2}_{L^\infty(\mathcal{Q}_{R,T})} \\  \ \\
\displaystyle \hspace{2cm}  +D^{\frac{1}{3}(\alpha-1)}  \|\nabla q\|^{1/3}_{L^\infty(\mathcal{Q}_{R,T})}\Big) 
  \Big(\beta + \ln \frac{D}{u(x,t)}\Big),
\end{array} \right.
\end{equation}
\item if $0< \alpha < 1$
\begin{equation}\label{eq12b}
\left. \begin{array}{l}
\displaystyle \frac{|\nabla u(x,t)|}{u(x,t)} \leq  C(\delta)  \Big( \sqrt{\frac{1+|\mu|}{R}} + \frac{1}{\sqrt{t-(t_0-T))}} +  \sqrt{K}   +  \sqrt{\alpha} \mathbb{M}^{\frac{1}{2}(\alpha-1)}\|q^+\|^{1/2}_{L^\infty(\mathcal{Q}_{R,T})} \\  \ \\
\displaystyle \hspace{2cm}  + \mathbb{M}^{\frac{1}{3}(\alpha-1)}  \|\nabla q\|^{1/3}_{L^\infty(\mathcal{Q}_{R,T})}\Big) 
  \Big(\beta + \ln \frac{D}{u(x,t)}\Big),
\end{array} \right.
\end{equation}
where $\mathbb{M}:= \inf \{u(x,t) : \text{for all} \  (x,t) \in \mathcal{Q}_{R, T}\}$,
\item if $\alpha \leq 0$
\begin{equation}\label{eq13}
\left. \begin{array}{l}
\displaystyle \frac{|\nabla u(x,t)|}{u(x,t)} \leq  C(\delta)  \Big( \sqrt{\frac{1+|\mu|}{R}} + \frac{1}{\sqrt{t-(t_0-T))}} +  \sqrt{K}   +  \mathbb{M}^{\frac{1}{2}(\alpha-1)}\|q^+\|^{1/2}_{L^\infty(\mathcal{Q}_{R,T})} \\  \ \\
\displaystyle \hspace{2cm}  + \mathbb{M}^{\frac{1}{3}(\alpha-1)}  \|\nabla q\|^{1/3}_{L^\infty(\mathcal{Q}_{R,T})}\Big) 
  \Big(\beta + \ln \frac{D}{u(x,t)}\Big).
\end{array} \right.
\end{equation}
 Here $q^+(x)= \max\{q(x),0\}$ and $\mu := \max_{\{x | d(x,x_0)=1\}} \Delta_f r(x)$, where $r(x)$ is the distance from a fixed point $x_0$ to point $x$ in $M$.
\end{enumerate}
\end{theorem}
 
Notice that when $q(x,t) \equiv 0$ (\ref{eq11}) reduces to the weighted heat equation
\begin{align}\label{eq11c}
\displaystyle \Big(\Delta_f - \frac{\partial}{\partial t}\Big) u(x,t) = 0.
\end{align} 
As an application of Theorem \ref{thm11}, we derive some Liouville-type theorems for positive ancient solutions to (\ref{eq11}) and (\ref{eq11c}) under certain  growth condition near infinity, when $Ric_f \geq 0$. This result is similar to the case on manifold with nonnegative Ricci tensor obtained by X. Zhu \cite{[Zhu16]}.

\begin{theorem}\label{thm12}
Let $(M^N, g, e^{-f}dv)$ be an $N$-dimensional complete smooth metric measure space with $Ric_f \geq 0$. Suppose that $q(x,t)=q(x)$, that is, time-independent and satisfies the following conditions
\begin{enumerate}[label=(\alph*)]
\item $\|q^+\|_{L^\infty(B(x_0,R))} = o(R^{-((\alpha-1)}) \ \ \ \text{as} \ \ \ R \to \infty$
\item $\|\nabla q\|_{L^\infty(B(x_0,R))} = o(R^{-(\alpha-1)}) \ \ \ \text{as} \ \ \ R \to \infty$.
\end{enumerate}
Then;
\begin{enumerate}[label=(\arabic*)]
\item For $q(x) \not\equiv 0$ equation (\ref{eq11}) has no positive ancient solution with $u(x,t) = o([r^{1/2}(x) +|t|^{1/4}])$ near infinity. (Ancient solution is a solution defined in all space and negative time).

Furthermore, for $q(x,t) \equiv 0$,  
\item equation (\ref{eq11c}) has only constant positive ancient solution with $u(x,t) = o([r^{1/2}(x) +|t|^{1/4}])$ near infinity,
\item equation (\ref{eq11c}) has only constant  ancient solution with $u(x,t) = o([r^{1/2}(x) +|t|^{1/4}])$ near infinity,

where $r$ is the distance from $x$ to a fixed point $y \in M$.
\end{enumerate}
\end{theorem}

Another application of our gradient estimates is the following Harnack type-inequalities:

\begin{theorem}\label{thm13}
Let $(M^N, g, e^{-f}dv)$ be an $N$-dimensional complete smooth metric measure space with $Ric_f \geq -(N-1)K$ for some $K \geq 0$.
Suppose  $\|q^+\|_{L^\infty}< \infty$ and $\|\nabla q\|_{L^\infty} < \infty$.
If $u(x,t)$ is a positive solution to (\ref{eq11}) and $u \leq D$ for all $(x,t) \in M \times [0, \infty),$ then 
\begin{align}
u(y,t) \leq u(x,t)^{\Gamma(r(x,y),t)}(De)^{1-\Gamma(r(x,y),t)}
\end{align}   
for all $x,y \in M$, where 
$$\Gamma(r(x,y),t)= \exp\Big(-C(\delta)\Big(\frac{1}{\sqrt{t-(t_0-T)}} + \sqrt{K} + \lambda \Big)r\Big),$$
where
$$\lambda :=  \max \{\sqrt{\alpha} D^{\frac{1}{2}(\alpha-1)}\|q^+\|^{1/2}_{L^\infty(\mathcal{Q}_{R,T})}, \ D^{\frac{1}{3}(\alpha-1)}  \|\nabla q\|^{1/3}_{L^\infty(\mathcal{Q}_{R,T})}\} \ \ \text{for} \ \ \alpha \geq 1,$$
$$\lambda :=  \max \{\sqrt{\alpha} \mathbb{M}^{\frac{1}{2}(\alpha-1)}\|q^+\|^{1/2}_{L^\infty(\mathcal{Q}_{R,T})}, \ \mathbb{M}^{\frac{1}{3}(\alpha-1)}  \|\nabla q\|^{1/3}_{L^\infty(\mathcal{Q}_{R,T})}\} \ \ \text{for} \ \ 0 < \alpha < 1,$$
$$\lambda :=  \max \{ \mathbb{M}^{\frac{1}{2}(\alpha-1)}\|q^+\|^{1/2}_{L^\infty(\mathcal{Q}_{R,T})}, \ \mathbb{M}^{\frac{1}{3}(\alpha-1)}  \|\nabla q\|^{1/3}_{L^\infty(\mathcal{Q}_{R,T})}\} \ \ \text{for} \ \ \alpha \leq 0 ,$$
$\mathbb{M} := \inf \{u(x,t) : \text{for all} \  (x,t) \in M \times [0, \infty) \}$ and  $r=r(x,y)$ denotes the geodesic distance between $x$ and $y$.
\end{theorem}
\begin{remark}
Note that we study (\ref{eq11}) for smooth function $f$ and an arbitrary constant $\alpha$. The elliptic gradient estimates obtained extend and generalize some known results, for instance,
 Theorem 1.1 of \cite{[Zhu11]}, where (\ref{eq11}) was studied for constant function $f$ and $\alpha \in (0,1)$, 
 Theorem 1.7 of \cite{[Zhu16]}, where (\ref{eq11}) was studied for constant function $f$ and $\alpha >1$ and 
 Theorem 1.2 of \cite{[MH17]}, where (\ref{eq11}) was studied for  arbitrary constants $q$ and $\alpha$.
\end{remark}
\begin{remark}
Puting $q(x,t)\equiv 0$, in (\ref{eq11}) we can deduce gradient estimate (1.3) of Theorem 1.1 obtained in \cite{[Wu15]} for the weighted heat equation. Obviously, our estimate (\ref{eq12}) (resp.(\ref{eq12b}) and (\ref{eq13})) generalize estimate (1.3) of Theorem 1.1  in \cite{[Wu15]}. Also, our estimates (\ref{eq12}) and (\ref{eq13}) generalize estimates (2.1) and (2.3) of Theorem 2.1 in \cite{[Wa15]} with $f$ being a constant and $\alpha \geq 1$ and $\alpha \leq 0$.
\end{remark}
The rest of this paper is organized as follows. In Section \ref{sec2}, we will give
some basic lemmas and introduce a space-time cut-off function that will be used in the proofs of Theorems \ref{thm11} and \ref{thm12}.
Section \ref{sec3} presents detail proofs of main results.

\section{Basic lemma}\label{sec2}
In this section we apply the arguments in Wu \cite{[Wu15]} and Zhu \cite{[Zhu16]} to prove Theorem \ref{thm11}.
Define a smooh function $h(x,t) = \ln u(x,t)/D$ for some constant $D$ in $\mathcal{Q}_{R,T}$. It is obvious that $h\leq 0$ and $h$ satisfies

\begin{align}\label{eq21}
\Big(\Delta_f - \frac{\partial}{\partial t}\Big) h(x,t) +|\nabla h(x,t)|^2 +q (De^{h})^{\alpha-1} =0.  
\end{align}
With this we prove the following lemma which is an extension of \cite{[SZ06],[Wu10b]}.

\begin{lemma}\label{lem21}
Let $(M^N, g, e^{-f}dv)$ be an $N$-dimensional smooth complete metric measure space with $Ric_f \geq -(N-1)K$ for some $K \geq 0$. Fix $x_0 \in M$ and $R \geq 0$. 
Let $h=h(x,t)$ be a smooth non-positive solution to (\ref{eq21})  in $\mathcal{Q}_{R,T}$, Then for all $(x,t) \in \mathcal{Q}_{R,T}$,  the function 
\begin{align}\label{eq22}
w = |\nabla \ln (\beta-h)|^2 = \frac{|\nabla h|^2}{(\beta-h)^2}
\end{align}
 satisfies
\begin{equation}\label{eq23}
\left. \begin{array}{l}
\displaystyle\Big(\Delta_f - \frac{\partial}{\partial t}\Big)w \geq  \frac{2[h+(1-\beta)]}{\beta -h}\langle \nabla h, \nabla w\rangle +2(\beta-h)w^2 - 2(N-1)Kw \\  \ \\
\displaystyle \hspace{2cm} - 2\Big[\alpha + \frac{h}{\beta -h} + \frac{1-\beta}{\beta -h}\Big] (De^{h})^{\alpha-1}q w - \frac{2}{(\beta -h)^2} (De^{h})^{\alpha-1}\langle \nabla h, \nabla q\rangle
\end{array} \right.
\end{equation}
for all $(x,t)$ in $\mathcal{Q}_{R,T}$,
\end{lemma}

\proof
We mostly work in a local orthonormal system with the convention that $h^2_i = |\nabla h|^2, h_{ii} = \Delta h$ and $h_{ijk}$ is the third order covariant derivative, and repeated indices are summed up.

Using (\ref{eq21}) and (\ref{eq22}) we compute that
\begin{align*}
\displaystyle w_t &  = \frac{2h_i(h_t)_i}{(\beta-h)^2} + \frac{2h^2_i h_t}{(\beta-h)^3}\\ 
\displaystyle &= \frac{2h_i(\Delta_f h+h_j^2+q (De^h)^{\alpha-1})_i}{(\beta-h)^2} + \frac{2h^2_i (\Delta_f h+h_j^2+q(De^h)^{\alpha-1})}{(\beta-h)^3}\\  
\displaystyle &= \frac{2h_i(\Delta_fh)_i}{(\beta-h)^2} +\frac{4h_ih_jh_{ij}}{(\beta-h)^2} + \frac{2h_iq_i(De^h)^{\alpha-1}}{(\beta-h)^2} + \frac{2(\alpha -1)(De^h)^{\alpha-1} h_i^2 q}{(\beta - h)^2} +  \frac{2h^2_i \Delta_fh}{(\beta-h)^3}\\
\displaystyle & \hspace{3cm} + \frac{2h^2_ih_j^2}{(\beta-h)^3} +  \frac{2h^2_i q (De^h)^{\alpha-1}}{(\beta-h)^3}.
\end{align*}

Similarly,
$$w_j  = \frac{2h_ih_{ij}}{(\beta-h)^2} + \frac{2h^2_i h_j}{(\beta-h)^3}$$
and
\begin{align*}
\displaystyle w_{jj} & =  \Big(\frac{2h_ih_{ij}}{(\beta-h)^2}\Big)_j  + \Big(\frac{2h^2_i h_j}{(\beta-h)^3}\Big)_j\\
\displaystyle & = \frac{2h^2_{ij}}{(\beta-h)^2} + \frac{2h_ih_{ijj}}{(\beta-h)^2} +\frac{8h_ih_jh_{ij}}{(\beta-h)^3}+\frac{2h^2_i h_{jj}}{(\beta-h)^3} + \frac{6h^2_i h^2_j}{(\beta-h)^4}.
\end{align*}
Using the following Ricci identity
$h_{ijj} = h_{jji}+R_{ij}h_j$,
a straight forward computation yields
$$\frac{2h_ih_{ijj}}{(\beta-h)^2} - \frac{2h_ih_{ij} f_j}{(\beta-h)^2} = \frac{2h_i (h_{jj} - h_jf_j)_i}{(\beta-h)^2} + \frac{2(R_{ij}+f_{ij})h_ih_j}{(\beta-h)^2}. $$
Hence
\begin{align*}
\displaystyle \Delta_f w & = \Delta w - \langle\nabla f, \nabla w \rangle  = w_{jj} - w_jf_j\\
\displaystyle & = \frac{2h^2_{ij}}{(\beta-h)^2} + \frac{2h_i (h_{jj} - h_jf_j)_i}{(\beta-h)^2}
+\frac{2(R_{ij}+f_{ij})h_ih_j}{(\beta-h)^2}
 +\frac{8h_ih_jh_{ij}}{(\beta-h)^3}\\
\displaystyle & \hspace{2cm}+ \frac{2h^2_i(h_{jj}-h_jf_j)}{(\beta-h)^3} + \frac{6h^2_i h^2_j}{(\beta-h)^4}\\
 \displaystyle & = \frac{2h^2_{ij}}{(\beta-h)^2} + \frac{2h_i(\Delta_f h)_i}{(\beta-h)^2}
+\frac{2 Ric_f(\nabla h, \nabla h)}{(\beta-h)^2}
 +\frac{8h_ih_jh_{ij}}{(\beta-h)^3}+\frac{2h^2_i\Delta_f h}{(\beta-h)^3} + \frac{6h^2_i h^2_j}{(\beta-h)^4}.
\end{align*}
Combining the above computations for $\Delta_fw$ and $w_t$, we have
\begin{align*}
\displaystyle \Big(\Delta_f-\frac{\partial}{\partial t}\Big) w  & = \frac{2h^2_{ij}}{(\beta-h)^2} 
+\frac{2(R_{ij}+f_{ij})f_if_j}{(\beta-h)^2}
 +\frac{8h_ih_jh_{ij}}{(\beta-h)^3} + \frac{6h^2_i h^2_j}{(\beta-h)^4}  - \frac{4h_ih_jh_{ij}}{(\beta-h)^2} \\
\displaystyle &   - \frac{2 (De^h)^{\alpha-1}h_iq_i}{(\beta-h)^2}
 - \frac{2h^2_ih_j^2}{(\beta-h)^3} - \frac{ 2 (De^h)^{\alpha-1}h^2_i q}{(\beta-h)^3}  - \frac{2(\alpha-1)h^2_i q (De^h)^{\alpha-1}}{(\beta-h)^3} \\
 & = \Big(\frac{2h^2_{ij}}{(\beta-h)^2} +\frac{4h_ih_jh_{ij}}{(\beta-h)^3} + \frac{2h^2_i h^2_j}{(\beta-h)^4}\Big) + \Big(\frac{4h_ih_jh_{ij}}{(\beta-h)^3} + \frac{4h^2_i h^2_j}{(\beta-h)^4} - \frac{4h_ih_jh_{ij}}{(\beta-h)^2} \\
 & - \frac{2h^2_i h^2_j}{(\beta-h)^3}\Big) +\frac{2(R_{ij}+f_{ij})f_if_j}{(\beta-h)^2}
- \frac{2(De^h)^{\alpha-1}h_iq_i}{(\beta-h)^2}
- \frac{2(De^h)^{\alpha-1}h^2_i q}{(\beta-h)^3}\\
&  - \frac{2(\alpha-1)h^2_i q (De^h)^{\alpha-1}}{(\beta-h)^3} \\  
& \geq \Big(\frac{2}{\beta-h}h_jwj - 2 h_jw_j + \frac{2h_i^2h_j^2}{(\beta-h)^3} \Big) - 2(N-1)K\frac{|\nabla h|^2}{(\beta-h)^2}\\
 & - \frac{2(De^h)^{\alpha-1}h_iq_i}{(\beta-h)^2}
- \frac{2(De^h)^{\alpha-1}h^2_i q}{(\beta-h)^3}  - \frac{2(\alpha-1)h^2_i q (De^h)^{\alpha-1}}{(\beta-h)^3} \\  
 & = \frac{2[h+(1-\beta)]}{\beta-h}\langle \nabla h, \nabla w\rangle + 2(\beta-h)w^2 - 2(N-1)Kw \\
 &  - 2(\alpha-1)(De^h)^{\alpha-1} q w -  \frac{2(De^h)^{\alpha-1}}{\beta-h} q w - \frac{2}{(\beta -h)^2} (De^{h})^{\alpha-1}\langle \nabla h, \nabla q \rangle  
\end{align*}
where we have used the condition $Ric_f \geq -(N-1)K$ and the following identities
$$\langle \nabla h, \nabla w\rangle = h_jw_j =  \frac{2h_ih_{ij}h_j}{(\beta-h)^2} + \frac{2h^2_i h^2_j}{(\beta-h)^3}$$
and
$$\Big(\frac{2h^2_{ij}}{(\beta-h)^2} +\frac{4h_ih_jh_{ij}}{(\beta-h)^3} + \frac{2h^2_i h^2_j}{(\beta-h)^4}\Big)= 2\Big(\frac{h_{ij}}{\beta-h} + \frac{h_ih_j}{(\beta-h)^2}\Big)^2 \geq 0. $$
This concludes the proof.

\qed

To prove Theorem \ref{thm11}, we shall apply the last Lemma and the localisation technique of Souplet-Zhang \cite{[SZ06]}. The theorem gives the elliptic gradient estimate for the positive smooth solutions to the linear weighted evolution equation (\ref{eq11}). We first introduce a well known cut-off function taken from \cite{[Wu15]} (see also \cite{[SZ06]}). The cut-off estimates will allow us derive the desired bounds in $\mathcal{Q}_{R,T}$.
\begin{lemma}\label{lem22}
Fix $t_0 \in \mathbb{R}$ and $T>0$. For any given $\tau \in (t_0-T, t_0]$, there exists a smooth function $\psi : [0, \infty)\times [t_0-T,t_0]\to \mathbb{R}$ satisfying the following properties
\begin{enumerate}[label=(\arabic*)]
\item $\psi = \psi(d(x,x_0),t) \equiv \psi(r, t); \ \psi(r,t) =1$ in $\mathcal{Q}_{R/2, T/2}$, \ $0\leq \psi(r,t) \leq 1$.
\item $\psi$ is a radially decreasing function in spatial variables and $\frac{d \psi}{dr}=0$ in $\mathcal{Q}_{R/2,T}.$ 
\item $|\frac{\partial \psi}{\partial r}|\frac{1}{\psi^a} \leq \frac{C_a}{R} \ \ \text{and} \ \ \ |\frac{\partial^2 \psi}{\partial r^2}|\frac{1}{\psi^a} \leq \frac{C_a}{R^2}$ in $[0, \infty)\times [t_0-T,t_0]$, where $0<a<1$.
\item $|\frac{\partial \psi}{\partial t}|\frac{1}{\psi^{1/2}} \leq \frac{C}{\tau-(t_0-T)}$ in $[0, \infty)\times [t_0-T,t_0]$ for some constant $C>0$ and $\psi(r, t_0-T)=0$ for all $r \in [0, \infty)$.
\end{enumerate}
\end{lemma}
We shall now apply Lemma \ref{lem21} and \ref{lem22} to prove Theorem \ref{thm11} via the maximum principle in a
local space-time supported set. We  mainly follow the arguments in \cite{[Wu15]}.

\section{Proof of Theorems \ref{thm11}, \ref{thm12}, and \ref{thm13}}\label{sec3}
\subsection*{\it Proof of Theorem \ref{thm11}} 
Choose a smooth function $\psi$ with support in $\mathcal{Q}_{R,T}$ and satisfies Lemma \ref{lem21}.  We then estimate $(\Delta_f-\partial_t)(\psi w)$ and analyse the result at a space-time point where the function $\psi w$ attains its maximum. 

A straightforward computation yields
\begin{align}\label{eq24}
\displaystyle\Big(\Delta_f -  \frac{\partial}{\partial t}\Big)(\psi w) = \psi \Big(\Delta_f -  \frac{\partial}{\partial t}\Big)w +  2\nabla w \nabla \psi + w \Big(\Delta_f -  \frac{\partial}{\partial t}\Big) \psi. 
\end{align}
Using Lemma \ref{lem21} in (\ref{eq24}) we have 
\begin{align}\label{eq25}
\left. \begin{array}{l}
\displaystyle\Big(\Delta_f -  \frac{\partial}{\partial t}\Big)(\psi w) -\Big(d + 2 \frac{\nabla \psi}{\psi} \Big) \nabla(\psi w) \\
\displaystyle \hspace{2cm}  \geq 2\psi(\beta-h)w^2 - (d \cdot \nabla \psi)w - 2 \frac{|\nabla \psi|^2}{\psi}w  - 2(N-1)K \psi w \\ 
\displaystyle \hspace{2cm} - \frac{2}{(\beta-h)^2}(De^h)^{\alpha-1}\psi\langle \nabla h, \nabla q\rangle  - 2 \Big[\alpha + \frac{h }{\beta-h} + \frac{1-\beta}{\beta-h}\Big](De^h)^{\alpha-1} q\psi w \\
\displaystyle \hspace{2cm} + w \Big(\Delta_f - \frac{\partial }{\partial t} \Big)\psi.
 \end{array} \right.
\end{align}
where we have used the identity
\begin{align*}
(d \cdot \nabla \psi)\psi + 2 \nabla \psi \nabla w = \Big(d + 2 \frac{\nabla \psi}{\psi} \Big) \nabla(\psi w) - (d \cdot \nabla \psi)w  - 2 \frac{|\nabla \psi|^2}{\psi} w
\end{align*}
and 
$$d:= \frac{2[h+(1-\beta)]}{\beta-h}\nabla h.$$
Suppose the space-time maximum of $\psi w$ is attained at the point $(x_1,t_1)$ in $\mathcal{Q}_{R,T}$. We can assume without loss of generality that $x_1$ is not in the cut locus of $M$, due to Calabi's argument \cite{[LY86]}. We also assume that $(\psi w)(x_1, t_1) >0$, otherwise $w(x,t) \leq 0$ and the result holds trivially. Then at the point $(x_1,t_1)$ (which is the maximal) we have 
$$\Delta_f(\psi w) \leq 0, \ \ \ (\psi w)_t \geq 0 \  \ \text{and} \ \ \nabla(\psi w) =0.$$
By the last estimates at $(x_1,t_1)$, (\ref{eq25}) can be simplified as 
\begin{align}\label{eq26}
\left. \begin{array}{l}
\displaystyle  2\psi(1-h)w^2 \leq (d \cdot \nabla \psi)w  + 2 \frac{|\nabla \psi|^2}{\psi} w + 2(N-1)K \psi w \\ 
\displaystyle \hspace{2cm} + \frac{2}{(\beta-h)^2}(De^h)^{\alpha-1}\psi\langle \nabla h, \nabla q\rangle  +  2 \Big[\alpha + \frac{h }{\beta-h} + \frac{1-\beta}{\beta-h}\Big](De^h)^{\alpha-1} q\psi w \\
\displaystyle \hspace{2cm} - w \Big(\Delta_f - \frac{\partial }{\partial t} \Big)\psi.
 \end{array} \right.
\end{align}

We now have two situations to consider; namely, if $x \not \in B(x_0,1)$ and if $x \in B(x_0,1)$. Firstly, we consider the situation if $x \not \in B(x_0,1)$.
To do the analysis, we need upper bounds  for each term on the right hand side (RHS) of (\ref{eq26}) at $(x_1,t_1)$. 

Let $C$ be a constant depending only on $N$, $C(\delta)$ a constant depending on $N$ and $\delta$ and their values vary from line to line. 
Closely following the arguments in \cite{[SZ06],[Wu10b]} with repeated use of Young's inequality  and the condition that $\beta -h \geq \delta >0$.

For the first term on the RHS of (\ref{eq26}):
$$(d\cdot \nabla \psi)w= \Big(\frac{2h}{\beta-h} \nabla h \nabla \psi\Big) w + \Big(\frac{2(1-\beta)}{\beta-h} \nabla h \nabla \psi\Big) w,$$
then
\begin{align}\label{eq27}
\left. \begin{array}{l}
\displaystyle \Big(\frac{2h}{\beta-h} \nabla h \nabla \psi\Big) w  \leq 2 |h| |\nabla \psi|w^{\frac{3}{2}} = 2 [\psi(\beta-h)w^2]^{\frac{3}{4}} \cdot \frac{|h||\nabla \psi|}{[\psi(\beta-h)]^{3/4}}\\
\displaystyle \hspace{2cm} \leq \psi(\beta-h)w^2 + C\Big( \frac{h}{\beta-h}\cdot \frac{|\nabla \psi|}{\psi^{3/4}}\Big)^4\\
\displaystyle \hspace{2cm} \leq  \psi(\beta-h)w^2 + \frac{C}{R^4}\frac{h^4}{(\beta-h)^3}.
 \end{array} \right.
\end{align}
and 
\begin{align}\label{eq28}
\left. \begin{array}{l}
\displaystyle \Big(\frac{2(1-\beta)}{\beta-h} \nabla h \nabla \psi\Big) w  \leq 2 |1-\beta| |\nabla \psi|w^{3/2} = (\psi w^2)^{3/4} \cdot \frac{|1-\beta||\nabla \psi|}{\psi^{3/4}}\\
\displaystyle \hspace{2cm} \leq \frac{\delta}{8} (\psi w^2) + C(\delta) \Big(\frac{|\nabla \psi|}{\psi^{3/4}}\Big)^4\\
\displaystyle \hspace{2cm} \leq \frac{\delta}{8} (\psi w^2) + C(\delta)\frac{1}{R^4}.
 \end{array} \right.
\end{align}
For the second term on the RHS of (\ref{eq26}):
\begin{align}\label{eq29}
\left. \begin{array}{l}
\displaystyle \frac{2|\nabla \psi|^2}{\psi}w = 2 \psi^{1/2}w \cdot \frac{|\nabla \psi|^2}{\psi^{3/2}}\\
\displaystyle \hspace{2cm} \leq  \frac{\delta}{8} (\psi^{1/2} w)^2 +  C(\delta)\Big(\frac{|\nabla \psi|^2}{\psi^{3/2}} \Big)^2\\
\displaystyle \hspace{2cm} \leq  \frac{\delta}{8} (\psi^{1/2} w)^2 + C(\delta)\frac{1}{R^4}.
 \end{array} \right.
\end{align}
For the third term on the RHS of (\ref{eq26}):
\begin{align}\label{eq210}
\left. \begin{array}{l}
\displaystyle 2(N-1)K \psi w  \leq \frac{\delta}{8} (\psi^{1/2} w)^2 + C(\delta)((N-1)K  \psi^{1/2})^2 \\
\displaystyle \hspace{4.8cm} \leq \frac{\delta}{8}\psi w^2 + C(\delta)K^2, 
 \end{array} \right.
\end{align}
For the sixth term on the RHS of (\ref{eq26}).
$$- w \Big(\Delta_f - \frac{\partial }{\partial t} \Big)\psi = -w(\Delta_f \psi) + w \psi_t.$$
 By the property that $\psi$ is a radial function and the weighted Laplacian comparison theorem \cite[Theorem 3.1]{[WeW09]}
$$\Delta_fr(x) \leq \mu +(N-1)K(R-1),$$
where $r(x,x_0)\geq 1$ in $B(x_0,R)$, $\mu := \max_{\{x|d(x,x_0)=1\}} \Delta r(x)$ and $Ric_f \geq - (N-1)K$:
\begin{align}\label{eq213}
\left. \begin{array}{l}
\displaystyle -(\Delta_f \psi)w = -\Bigg(\Big(\frac{\partial \psi}{\partial r}\Big)\Delta_f r + \Big(\frac{\partial^2\psi}{\partial r^2}\Big)\cdot |\nabla r|^2 \Bigg)w \\
\displaystyle \hspace{2cm} \leq - \Bigg(\Big(\frac{\partial \psi}{\partial r}\Big)(\mu +(N-1)K(R-1)) + \frac{\partial^2\psi}{\partial r^2} \Bigg)w\\
\displaystyle \hspace{2cm} \leq \Bigg(\Big|\frac{\partial^2\psi}{\partial r^2}\Big| \frac{1}{\psi^{1/2}} + (\mu +(N-1)K(R-1))\Big|\frac{\partial \psi}{\partial r}\Big| \frac{1}{\psi^{1/2}}  \Bigg)\psi^{1/2}w\\
\displaystyle \hspace{2cm} \leq  \frac{\delta}{8} \psi w^2 + C(\delta)\Bigg(\Big|\frac{\partial^2\psi}{\partial r^2}\Big| \frac{1}{\psi^{1/2}} + (\mu +(N-1)K(R-1))\Big|\frac{\partial \psi}{\partial r}\Big| \frac{1}{\psi^{1/2}}  \Bigg)\\
\displaystyle \hspace{2cm} \leq  \frac{\delta}{8} \psi w^2 +\frac{C(\delta)}{R^4} + \frac{C(\delta)\mu^2}{R^2} + C(\delta)K^2
 \end{array} \right.
\end{align}
where we have used the property (4) in Lemma \ref{eq22},
and
\begin{align}\label{eq214}
\left. \begin{array}{l}
\displaystyle \psi_t w \leq \psi^{1/2}w\frac{|\psi_t|}{\psi^{1/2}}\\
\displaystyle \hspace{1cm} \leq  \frac{\delta}{8}\psi w^2 + C(\delta) \Big(\frac{|\psi_t|}{\psi^{1/2}}\Big)^2\\
\displaystyle \hspace{1cm} \leq \frac{1}{8}\psi w^2  + \frac{C(\delta)}{(\tau-(t_0-T))^2}.
 \end{array} \right.
\end{align}

{\bf Case 1: For $\alpha \geq 1$.}
  
For the fourth term on the RHS of (\ref{eq26}): Here we know that $0 < e^{h(\alpha-1)} \leq 1$ since $h$ is nonpositive. Therefore
\begin{align}\label{eq211}
\left. \begin{array}{l}
\displaystyle \frac{2}{(\beta-h)^2} (De^h)^{(\alpha-1)} \psi\langle \nabla h, \nabla q \rangle \leq   \frac{2}{(\beta-h)^2} (De^h)^{(\alpha-1)} \psi  |\nabla h|| \nabla q |\\
\displaystyle \hspace{5.5cm} \leq \frac{2}{(\beta-h)} D^{\alpha-1} \psi | \nabla q | w^{1/2}\\
\displaystyle \hspace{5.5cm} \leq \frac{\delta}{8}(\psi^{1/4} w^{1/2})^4 + C(\delta) \Big(\psi^{3/4}D^{\alpha-1} \frac{|\nabla q|}{\beta-h} \Big)^{4/3}\\
\displaystyle \hspace{5.5cm} \leq \frac{\delta}{8}\psi w^2 + C(\delta) D^{\frac{4}{3}(\alpha-1)} \frac{|\nabla q|^{4/3}}{(\beta-h)^{4/3}}\\
\displaystyle \hspace{5.5cm} \leq \frac{\delta}{8} \psi w^2 + C(\delta) D^{{\frac{4}{3}}(\alpha-1)} |\nabla q|^{4/3}.
 \end{array} \right.
\end{align}

For the fifth term on the RHS of (\ref{eq26}): Considering the conditions    $\beta-h \geq \delta >0$, $\beta := \sup_{(x,t) \in \mathcal{Q}_{R/2,T}}|h| +1 \geq 1$ and the fact that $h$ is nonpositive, we note that  
$$0 < e^{h(\alpha-1)} \leq 1, \ \ \ 0< \frac{-h}{\beta-h}  = 1 - \frac{\beta}{\beta-h} <1
\ \ \ \text{and} \ \ \ \frac{1-\beta}{\beta-h} \leq 0.$$
Therefore  
$$ -1 < \frac{h}{\beta-h}< 0 \ \ \ \text{and} \ \ \ 0< \alpha + \frac{h}{\beta-h} + \frac{1-\beta}{\beta-h}   <   \alpha
$$
and
\begin{align}\label{eq212}
\left. \begin{array}{l}
\displaystyle 2\Big[\alpha+ \frac{h}{\beta-h} + \frac{1-\beta}{\beta-h}\Big] (De^h)^{(\alpha-1)} q \psi w  \leq 2\Big[ \alpha + \frac{h}{\beta-h}  + \frac{1-\beta}{\beta-h}\Big] D^{(\alpha-1)} q^+ \psi w\\
\displaystyle \hspace{3.5cm}  \leq \frac{\delta}{8} \psi w^2 + C(\delta) \Big[\alpha + \frac{h}{\beta-h} + \frac{|1-\beta|}{\beta-h}\Big]^2 \Big(D^{(\alpha-1)} \psi^{1/2} q^+ \Big)^2\\
\displaystyle \hspace{3.5cm} \leq \frac{\delta}{8} \psi w^2 + C(\delta)\alpha^2 D^{2(\alpha-1)} (q^+)^2, 
\end{array} \right.
\end{align}
where $q^+(x) = \max\{q(x),0\}$.

Now substituting (\ref{eq27})--(\ref{eq212}   
 into the RHS of (\ref{eq26}) and rearranging we obtain
\begin{align}\label{eq215}
\left. \begin{array}{l}
\displaystyle  \psi(\beta-h)w^2 \leq  \frac{7 \delta}{8} \psi w^2 + \frac{C}{R^4}\frac{h^4}{(\beta-h)^4}  + \frac{C(\delta)}{R^4}   + C(\delta) D^{4/3(\alpha-1)}|\nabla q|^{4/3} 
\\
\displaystyle \hspace{2cm} + C(\delta) \alpha^2 D^{2(\alpha-1)}(q^+)  + \frac{C(\delta) \mu^2}{R^2} + C(\delta)K^2 + \frac{C(\delta)}{(\tau-(t_0-T))^2}   
 \end{array} \right.
\end{align}
at $(x_1,t_1)$.  
Since  $\beta-h \geq \delta >0$ and $h /(\beta - h) \in (-1,0)$   implies  $h^4/(\beta-h)^4 < 1,$ then (\ref{eq215}) implies

\begin{align*}
\left. \begin{array}{l}
\displaystyle (\psi w^2)(x_1,t_1) \leq  \frac{7}{8} \psi w^2 + \frac{C}{R^4} + C(\delta) D^{3/4(\alpha-1)}|\nabla q|^{4/3}  + C(\delta) \alpha^2 D^{2(\alpha-1)}(q^+) \\
\displaystyle \hspace{3.5cm} + \frac{C(\delta) \mu^2}{R^2} + C(\delta)K^2 + \frac{C(\delta)}n{(\tau-(t_0-T))^2}   \\
\displaystyle \hspace{2.7cm} \leq C(\delta)\Big( \frac{1}{R^4} + \frac{\mu^2}{R^2}  + K^2 + \frac{1}{(\tau-(t_0-T))^2}   +\alpha^2  D^{2(\alpha-1)} \|q^+\|^2_{L^\infty(\mathcal{Q}_{R,T})} \\  
\displaystyle \hspace{3.5cm} + D^{4/3(\alpha-1)} \|\nabla q\|^{4/3}_{L^\infty(\mathcal{Q}_{R,T})}\Big).  
 \end{array} \right.
\end{align*}
at $(x_1, t_1).$
It follows that for all $(x,t) \in \mathcal{Q}_{R,T}$, there holds
\begin{align}\label{eq216}
\left. \begin{array}{l}
\displaystyle  (\psi^2 w^2)(x,\tau) \leq (\psi^2 w^2)(x_1,t_1) \leq (\psi w^2)(x_1,t_1) \\
\displaystyle \hspace{2.5cm} \leq C(\delta)\Big( \frac{1}{R^4} + \frac{\mu^2}{R^2}  + K^2 + \frac{1}{(\tau-(t_0-T))^2}   +\alpha^2  D^{2(\alpha-1)} \|q^+\|^2_{L^\infty(\mathcal{Q}_{R,T})} \\ \ \\  
\displaystyle \hspace{3.5cm} + D^{4/3(\alpha-1)} \|\nabla q\|^{4/3}_{L^\infty(\mathcal{Q}_{R,T})}\Big).   
 \end{array} \right.
\end{align}
Note that  $\psi(x,\tau)=1$ in $\mathcal{Q}_{R/2, T/2}$ (by Lemma \ref{lem22}), $w = |\nabla h|^2/(\beta-h)^2$ (by definition) and the fact that $\tau \in (t_0-T, t_0]$ was arbitrarily chosen, we have 
\begin{align}\label{eq217}
\left. \begin{array}{l}
\displaystyle \frac{|\nabla h|}{(\beta-h)}(x,t)  \leq C(\delta)\Big( \frac{1}{R^4} + \frac{\mu^2}{R^2}  + K^2 + \frac{1}{(\tau-(t_0-T))^2}   +\alpha^2  D^{2(\alpha-1)} \|q^+\|^2_{L^\infty(\mathcal{Q}_{R,T})} \\  
\displaystyle \hspace{5cm} + D^{4/43\alpha-1)} \|\nabla q\|^{4/3}_{L^\infty(\mathcal{Q}_{R,T})}\Big)^{1/4} \\
\displaystyle \hspace{2.5cm}   \leq C(\delta)\Big( \frac{1}{R} + \sqrt{\frac{\mu}{R}}  + \sqrt{K} + \frac{1}{\sqrt{(t-(t_0-T))}}   + \sqrt{\alpha}  D^{1/2(\alpha-1)} \|q^+\|^{1/2}_{L^\infty(\mathcal{Q}_{R,T})} \\   
\displaystyle \hspace{5cm} + D^{1/3(\alpha-1)} \|\nabla q\|^{1/3}_{L^\infty(\mathcal{Q}_{R,T})}\Big)   
 \end{array} \right.
\end{align}
for all $(x,t) \in \mathcal{Q}_{R/2,T} \equiv B(x_0, R/2) \times[t_0-T, t_0]$ with $t \neq t_0-T$.  
Since $h =\ln u/D$, we have 
\begin{align}\label{eq218}
\displaystyle \frac{|\nabla h|}{(\beta-h)}(x,t) = \Big(\frac{|\nabla u|}{u}\frac{1}{\beta- \ln u/D} \Big)(x.t).
 \end{align} 
By substituting (\ref{eq218}) into (\ref{eq217}) and rearranging we arrive at (\ref{eq12}).

{\bf Case 2: For $0< \alpha < 1$ }. In this case We have $e^{h(\alpha-1)}>1$  since $h$ is nonpositive.
For the fourth term on the RHS of (\ref{eq26}):
\begin{align}\label{eq211k}
\left. \begin{array}{l}
\displaystyle \frac{2}{(\beta-h)^2} (De^h)^{(\alpha-1)} \psi\langle \nabla h, \nabla q \rangle \leq   \frac{2}{(\beta-h)^2} \mathbb{M}^{(\alpha-1)} \psi  |\nabla h|| \nabla q |\\
\displaystyle \hspace{5.5cm} \leq \frac{\delta}{8}(\psi^{1/4} w^{1/2})^4 + C(\delta) \Big(\psi^{3/4}\mathbb{M}^{\alpha-1} \frac{|\nabla q|}{\beta-h} \Big)^{4/3}\\
\displaystyle \hspace{5.5cm} \leq \frac{\delta}{8}\psi w^2 + C(\delta) \mathbb{M}^{\frac{4}{3}(\alpha-1)} \frac{|\nabla q|^{4/3}}{(\beta-h)^{4/3}}\\
\displaystyle \hspace{5.5cm} \leq \frac{\delta}{8} \psi w^2 + C(\delta) \mathbb{M}^{{\frac{4}{3}}(\alpha-1)} |\nabla q|^{4/3},
 \end{array} \right.
\end{align}

For the fifth term on the RHS of (\ref{eq26}): We have $e^{h(\alpha-1)}>1$  since $h$ is nonpositive. Since $\alpha \in (0,1)$ (and similarly to the case $\alpha \geq 1$), we have  
$$ \frac{h}{\beta-h}  \in (-1,0), \ \ \  \frac{1-\beta}{\beta-h} \leq 0 \ \ \ \text{and} \ \ \ 0 < \alpha+ \frac{h}{\beta-h} + \frac{1-\beta}{\beta-h}  \leq  \alpha.$$
Then
\begin{align}\label{eq212m}
\left. \begin{array}{l}
\displaystyle2\Big[\alpha+ \frac{h}{\beta-h} + \frac{1-\beta}{\beta-h}\Big] (De^h)^{(\alpha-1)} q \psi w  \leq 2\Big[ \alpha + \frac{h}{\beta-h}  + \frac{1-\beta}{\beta-h}\Big] \mathbb{M}^{(\alpha-1)} q^+ \psi w\\
\displaystyle \hspace{5.5cm}  \leq \frac{\delta}{8} \psi w^2 + C(\delta) \Big[\alpha + \frac{h}{\beta-h} + \frac{|1-\beta|}{\beta-h}\Big]^2 \Big(\mathbb{M}^{(\alpha-1)} \psi^{1/2} q^+ \Big)^2\\
\displaystyle \hspace{5.5cm}  \leq \frac{\delta}{8} \psi w^2 + C(\delta) \alpha^2 \Big(\mathbb{M}^{(\alpha-1)} \psi^{1/2} q^+ \Big)^2\\
\displaystyle \hspace{5.5cm} \leq \frac{\delta}{8} \psi w^2 + C(\delta)\alpha^2 \mathbb{M}^{2(\alpha-1)} (q^+)^2, 
\end{array} \right.
\end{align}
where $q^+(x) = \max\{q(x),0\}$.

Similarly, putting (\ref{eq27})--(\ref{eq214} and    (\ref{eq211k})--(\ref{eq212m}
 into the RHS of (\ref{eq26}) and rearranging we obtain 
\begin{align*}
\left. \begin{array}{l}
\displaystyle  \psi(\beta-h)w^2 \leq  \frac{7 \delta}{8} \psi w^2 + \frac{C}{R^4}\frac{h^4}{(\beta-h)^4}  + \frac{C(\delta)}{R^4}   + C(\delta) \mathbb{M}^{4/3(\alpha-1)}|\nabla q|^{4/3} 
\\
\displaystyle \hspace{2cm} + C(\delta) \alpha^2 \mathbb{M}^{2(\alpha-1)}(q^+)  + \frac{C(\delta) \mu^2}{R^2} + C(\delta)K^2 + \frac{C(\delta)}{(\tau-(t_0-T))^2}   
 \end{array} \right.
\end{align*}
at $(x_1,t_1)$.  It then follows that for all $(x,t) \in \mathcal{Q}_{R,T}$, there holds
\begin{align}
\left. \begin{array}{l}
\displaystyle  (\psi^2 w^2)(x,\tau) \leq (\psi^2 w^2)(x_1,t_1) \leq (\psi w^2)(x_1,t_1) \\
\displaystyle \hspace{2.5cm} \leq C(\delta)\Big( \frac{1}{R^4} + \frac{\mu^2}{R^2}  + K^2 + \frac{1}{(\tau-(t_0-T))^2}   +\alpha^2  \mathbb{M}^{2(\alpha-1)} \|q^+\|^2_{L^\infty(\mathcal{Q}_{R,T})} \\  
\displaystyle \hspace{3.5cm} + \mathbb{M}^{4/3(\alpha-1)} \|\nabla q\|^{4/3}_{L^\infty(\mathcal{Q}_{R,T})}\Big).   
 \end{array} \right.
\end{align}
Following similar steps as before we arrive at (\ref{eq12b}).

{\bf Case 3: For $\alpha \leq 0$}

For the fourth term on the RHS of (\ref{eq26}):
We know that $e^{h(\alpha-1)}>1$  since $h$ is nonpositive and 
\begin{align}\label{eq211l}
\left. \begin{array}{l}
\displaystyle \frac{2}{(\beta-h)^2} (De^h)^{(\alpha-1)} \psi\langle \nabla h, \nabla q \rangle \leq   \frac{2}{(\beta-h)^2} \mathbb{M}^{(\alpha-1)} \psi  |\nabla h|| \nabla q |\\
\displaystyle \hspace{5.5cm} \leq \frac{\delta}{8}(\psi^{1/4} w^{1/2})^4 + C(\delta) \Big(\psi^{3/4}\mathbb{M}^{\alpha-1} \frac{|\nabla q|}{\beta-h} \Big)^{4/3}\\
\displaystyle \hspace{5.5cm} \leq \frac{\delta}{8}\psi w^2 + C(\delta) \mathbb{M}^{\frac{4}{3}(\alpha-1)} \frac{|\nabla q|^{4/3}}{(\beta-h)^{4/3}}\\
\displaystyle \hspace{5.5cm} \leq \frac{\delta}{8} \psi w^2 + C(\delta) \mathbb{M}^{{\frac{4}{3}}(\alpha-1)} |\nabla q|^{4/3},
 \end{array} \right.
\end{align}
where $\mathbb{M}:= \inf \{u(x,t) : \text{for all} \  (x,t) \in \mathcal{Q}_{R, T}\}$.

For the fifth term on the RHS of (\ref{eq26}): we know also that $e^{h(\alpha-1)}>1$. By the condition  $\beta-h \geq \delta >0$ we have  $ h \leq \beta-\delta$, 
$$\frac{h}{\beta-h} \leq \frac{\beta}{\delta} -1, \ \ \ \ \frac{1-\beta}{\beta-h} \leq \frac{1}{\delta} - \frac{\beta}{\delta} \ \ \ \text{and} \ \ \ \alpha+ \frac{h}{\beta-h} + \frac{1-\beta}{\beta-h}   \leq   \frac{1}{\delta}-1 <  \frac{1}{\delta}.$$
Therefore
\begin{align}\label{eq212n}
\left. \begin{array}{l}
\displaystyle 2\Big[\alpha+ \frac{h}{\beta-h} + \frac{1-\beta}{\beta-h}\Big] (De^h)^{(\alpha-1)} q \psi w  \leq 2\Big[ \alpha + \frac{h}{\beta-h}  + \frac{1-\beta}{\beta-h}\Big] \mathbb{M}^{(\alpha-1)} q^+ \psi w\\
\displaystyle \hspace{3.5cm}  \leq \frac{\delta}{8} \psi w^2 + C(\delta) \Big[\alpha + \frac{h}{\beta-h} + \frac{|1-\beta|}{\beta-h}\Big]^2 \Big(\mathbb{M}^{(\alpha-1)} \psi^{1/2} q^+ \Big)^2\\
\displaystyle \hspace{3.5cm}  \leq \frac{\delta}{8} \psi w^2 + C(\delta) \Big(\mathbb{M}^{(\alpha-1)} \psi^{1/2} q^+ \Big)^2\\
\displaystyle \hspace{3.5cm} \leq \frac{\delta}{8} \psi w^2 + C(\delta) \mathbb{M}^{2(\alpha-1)} (q^+)^2, 
\end{array} \right.
\end{align}
where $q^+(x) = \max\{q(x),0\}$ and $\mathbb{M}:= \inf \{u(x,t) : \text{for all} \  (x,t) \in \mathcal{Q}_{R, T}\}$.
Similarly putting (\ref{eq27})--(\ref{eq214} and    (\ref{eq211l})--(\ref{eq212n}
 into the RHS of (\ref{eq26}), rearranging and following the same steps as before we obtain (\ref{eq13}).

Now we consider the other situation: if $x \in B(x_0,1)$. Here $\psi$ is a constant in space direction in $B(x_0,R/2)$ based on the assumption, where $R \geq 2$. Thus at $(x_1,t_1)$, we have from $(\ref{eq26})$ for the case $\alpha \geq 1$ (Note that $\beta-h\geq \delta$ and $0 < e^{h(\alpha-1)}\leq 1$)
\begin{align*}
\displaystyle w &\leq \frac{\psi_t}{2\psi} + (N-1)K + \frac{2}{\beta-h}(De^h)^{\alpha-1} |\nabla q|  +  2 \Big[\alpha + \frac{h }{\beta-h} + \frac{1-\beta}{\beta-h}\Big](De^h)^{\alpha-1} (q^+) \\
\displaystyle & \leq \frac{C}{\tau-(t_0-T)} + (N-1)K + C(\delta) D^{\alpha-1} |\nabla q|+ \alpha   D^{\alpha-1}(q^+),
\end{align*}
where we have used $(4)$ of Lemma \ref{lem22}. Since $\psi(x,\tau) =1$ when $d(x,x_0) < R/2$ by $(1)$ of Lemma \ref{lem22}, the  last estimate indeed yields 
\begin{align*}
\displaystyle w(x,\tau) = (\psi w)(x,\tau) & \leq (\psi w)(x_1,t_1)\\
\displaystyle &\leq  w(x_1,t_1)\\
\displaystyle & \leq \frac{C}{\tau-(t_0-T)} + (N-1)K + C(\delta) D^{\alpha-1} |\nabla q|+ \alpha   D^{\alpha-1}(q^+)
\end{align*}
for all $(x,t) \in \mathcal{Q}_{R/2,T}$ with $t\neq t_0-T$. This proves the estimate (\ref{eq12}).
Similarly for the cases $\alpha \in (0,1)$ and $ \alpha \leq 0$, we can easily obtain estimates (\ref{eq12b}) and (\ref{eq13}) respectively. This concludes the proof of Theorem \ref{thm11}.
\qed
 
\subsection*{\it Proof of Theorem \ref{thm12}}
To prove $(1)$ of the Theorem, we only consider the case $\alpha \geq 1$ since the case $\alpha \leq 1$ is similar.
Let 0$<u(x,t) \leq D$ be a positive ancient solution of (\ref{eq11}) with $ u(x,t) = o([r(x)^{1/2} + |t|^{1/4}])$ near infinity. Fix any point $(x_0,t_0)$ in space-time and
let $D_R:= \sup_{\mathcal{Q}_{\sqrt{R},R}}|u|$. Considering the function $U= u + 2D_{2R}$, we have $D_{2R} \leq U(x,t) \leq 3D_{2R}$, for ever $(x,t) \in \mathcal{Q}_{2\sqrt{R},4R}$.
Then  using Theorem \ref{thm11} for $U$ in the set $B(x_0,R)\times [t_0-R^2,t_0]$, we have 
\begin{align*}
\displaystyle \frac{|\nabla u(x_0,t_0)|}{u(x_0,t_0) +  2D_{2R}} \leq & \frac{C(\delta,|\mu|)}{R}  \Big( R^{\frac{1}{2}}\Big) +C(\delta,\alpha)o(R^{\frac{1}{2}(\alpha-1)})o(R^{-\frac{1}{2}(\alpha-1)}) + \\ & C(\delta)o(R^{\frac{1}{3}(\alpha-1)}) o(R^{-\frac{1}{2}(\alpha-1)})
\end{align*}
near infinity. Since $ D_{2R} = o(R)$ by asumption,  it follows that $|\nabla u(x_0,t_0)|=0$ by letting $R\to \infty$. Since $(x_0,t_0)$ is arbitrary, then $\nabla u(x,t)\equiv 0$ and $u$ must be  constant in space, i.e., $u(x,t)=u(t)$. Furthermore, from equation (\ref{eq11}), we have $q(x) = \widetilde{q}$ (a constant) and we obtain
\begin{align}\label{eq217m}
\frac{du(t)}{dt} = \widetilde{q} u^\alpha(t)
\end{align}
Integrating (\ref{eq217}) in the interval $(t,0]$ with $t<0$ we obtain
\begin{align}\label{eq218n}
u^{1-\alpha}(t) = u^{1-\alpha}(0) + (1-\alpha)\widetilde{q} t.
\end{align}
We can then prove that $\widetilde{q} =0$ which is a contradiction to the assumption that $q(x) \neq 0$.
By the hypothesis $(a)$ of the Theorem we have $\widetilde{q} \leq 0$. Letting $t \to -\infty$, we have $u^{1-\alpha}(t) <0$ (for $\alpha \geq 1$), which is impossible because $u$ is a positive solution. Therefore $\widetilde{q}=0$ is a contradition to $u = o([r^{1/2}(x)+|t|^{1/4}])$ near infinity. 

We prove $(2)$ and $(3)$ of the Theorem. For $q(x) \equiv 0$, (\ref{eq11}) becomes the weighted heat equation (\ref{eq11c}. Let $u(x,t)$ be an ancient solution of  (\ref{eq11c} with $u = o([r^{1/2}(x)+|t|^{1/4}])$ near infinity. Then from Theorem \ref{thm11}, we have $\nabla u \equiv 0$ so $u(x,t) = u(t)$. Then from (\ref{eq11c} we have 
$$\frac{du(t)}{dt} =0$$
which implies that $u$ is a constant. This completes the proof of Theorem \ref{thm13}.

\subsection*{\it Proof of Theorem \ref{thm13}}
Let $\gamma(s)$, $\gamma: [0,1] \to M$, be the minimal geodesic connecting $x$ and $y$  such that $\gamma(0) = x$ and $\gamma(1) = y$. Let $h =\ln u/D$. By letting $R \to \infty$ in Theorem \ref{thm11} we have
\begin{align*}
\frac{|\nabla u|}{u(\beta- \ln u/D)} \leq C(\delta)\Big(\frac{1}{\sqrt{t-(t_0-T)}} + \sqrt{K} + \lambda \Big),
\end{align*}
where
$$\lambda :=  \max \{\sqrt{\alpha} D^{\frac{1}{2}(\alpha-1)}\|q^+\|^{1/2}_{L^\infty(\mathcal{Q}_{R,T})}, \ D^{\frac{1}{3}(\alpha-1)}  \|\nabla q\|^{1/3}_{L^\infty(\mathcal{Q}_{R,T})}\} \ \ \text{for} \ \ \alpha \geq 1,$$
$$\lambda :=  \max \{\sqrt{\alpha} \mathbb{M}^{\frac{1}{2}(\alpha-1)}\|q^+\|^{1/2}_{L^\infty(\mathcal{Q}_{R,T})}, \ \mathbb{M}^{\frac{1}{3}(\alpha-1)}  \|\nabla q\|^{1/3}_{L^\infty(\mathcal{Q}_{R,T})}\} \ \ \text{for} \ \ 0 < \alpha < 1,$$
and
$$\lambda :=  \max \{ \mathbb{M}^{\frac{1}{2}(\alpha-1)}\|q^+\|^{1/2}_{L^\infty(\mathcal{Q}_{R,T})}, \ \mathbb{M}^{\frac{1}{3}(\alpha-1)}  \|\nabla q\|^{1/3}_{L^\infty(\mathcal{Q}_{R,T})}\} \ \ \text{for} \ \ \alpha \leq 0.$$

We now compute
\begin{align*}
\ln \frac{\beta-h(x,t)}{\beta-h(y,t)} &= \int_0^1 \frac{d\ln(\beta-h(\gamma(s),t))}{ds} ds\\
& \leq \int_0^1 |\dot{\gamma}| \frac{|\nabla u|}{u(\beta-\ln u/D)} ds\\
& \leq C(\delta) \Big(\frac{1}{\sqrt{t-(t_0-T)}} + \sqrt{K} + \lambda \Big)r. 
\end{align*}
Denote by
$$\Gamma = \Gamma(r(x,y),t) := \exp\Big(-C(\delta)\Big(\frac{1}{\sqrt{t-(t_0-T)}} + \sqrt{K} + \lambda \Big)r\Big)$$
the above inequality implies
$$\frac{\beta-h(x,t)}{\beta-h(y,t)} \leq \frac{1}{\Gamma}.$$
Hence, with some straightforward computation, we obtain
\begin{align*}
u(y,t) \leq u(x,t)^{\Gamma}\cdot (De)^{1-\Gamma}
\end{align*} 
which concludes the proof.
\qed


\section*{Acknowledgement}
The author  wishes to thank the anonymous referees for their useful comments. He also thanks Dr Jia-Yong Wu for helpful discussions he had with him at various stages of this  work.


\begin{thebibliography}{00}
\bibitem{[BCP10]}{M. B\v{a}ile\c{s}teanu, X. D. Cao, A.  Pulemotov}
\emph{Gradient estimates for the heat equation under Ricci flow}. 
J. Funct. Anal, 258 (2010), 3517--3542.
 
\bibitem{[BQ00]}{D. Bakry, Z. M. Qian},
 \emph{Some new results on eigenvectors via dimension, diameter and Ricci curvature}.
Adv. Math. 155,(2000), 98--153.

\bibitem{[BRS]} {L. Brandolini, M. Rigoli, A. G. Setti}, \emph{Positive solutions of Yamabe type equations on complete
manifolds and applications}, J. Funct. Anal. 160 (1998), 176--222.

\bibitem{[Br13]}{K. Brighton},
\emph{A Liouville-type theorem for smooth metric measure spaces}.
 J. Geom. Anal. 23 (2013), 562--570.
 
 \bibitem{[CaH]}{H. D. Cao}, \emph{Recent progress on Ricci solitons, Recent advances in geometric analysis}. 
Adv. Lect. Math. (ALM), Int. Press, Somerville, MA (2010), 11, 1--38.

\bibitem{[Cas1]}{J. S. Case}, {\it Conformal invariants measuring the best constants for Gagliardo-
Nirenberg-Sobolev inequalities}, Calc. Var. Partial Diff. Eqns. 48(2013), 507--526. 

\bibitem{[Cas2]}{J. S. Case}, {\it A Yamabe-type problem on smooth metric measure spaces},
J. Diff. Geom. 101 (2015), 467--505.

\bibitem{[CC09]}{L. Chen, W.Y. Chen}, 
		{\it Gradient estimates for a nonlinear parabolic equation on complete noncompact Riemannian manifolds}. Ann Global Anal Geom, 35(4) (2009), 397--404.
			
			 
\bibitem{[CY75]} {S. Y. Cheng, S. T, Yau},
\emph{Differential equations on Riemannian manifolds and their geometric applications}.
Commun. Pure Appl. Math. 28 (1975), 333--354.

\bibitem{[DD]}{M. Del Pino, J. Dolbeault}, {\it Best constants for Gagliardo-Nirenberg inequalities
and applications to nonlinear diffusions}, J. Math. Pures Appl. 81(2002),
847--875.  

\bibitem{[Ha93]}{R. Hamilton},
\emph{A matrix Harnack estimate for the heat equation}. 
Comm. Anal. Geom. 1 (1993), 113--126.

\bibitem{[Ha95]}{R. Hamilton},
\emph{The formation of singularities in the Ricci flow}.
Surv. Differ. Geom. 2,(1995), 7--136.

\bibitem{[Le]}{M. Ledoux}, { \it The geometry of Markov diffusion generators}. 
Ann.  fac.  sci, de Toulouse (6), 2(2000), 305--366

\bibitem{[Li05]}{X. D. Li},  
\emph{Liouville theorems for symmetric diffusion operators on complete Riemannian manifolds}.
 J. Math. Pure. Appl. 84 (2005), 1295--1361
 
\bibitem{[LY86]}{ P. Li, S-T. Yau},  
\emph {On the parabolic kernel of the Schr\"odinger operator}, 
Acta Math. 156 (1986), 353--364.

\bibitem{[MH17]}{ B. Ma, G. Huang},  
\emph {Hamilton-Souplet-Zhang's estimates for two weightd nonlinear parabolic equations}, 
Appl. Math. J. Chinese Univ. 32(3) (2017), 153--201.

\bibitem{[MRS]}{P. Mastrolia, M. Rigoli, A. G. Setti}, {\it Yamabe-type Equations on Complete, Noncompact Manifolds}, Progress in Mathematics 302, Birkhäuser Verlag, Basel, 2012. 
 
\bibitem{[Ru10]}{Q. H. Ruan},
\emph{Elliptic-type gradient estimate for Schr\"odinger equations on noncompact manifolds}. Bull. Lond. Math. Soc., 39(6) (2007), 982--988.

\bibitem{[SZ06]}{P. Souplet, Q. S, Zhang},
\emph{Sharp gradient estimate and Yau's Liouville theorem for the heat equation on noncompact manifolds}. 
Bull. London Math. Soc. 38 (2006), 1045--1053

\bibitem{[Wa15]}{W. Wang},
\emph{Complement of gradient estimates and
Liouville theorems for nonlinear parabolic
equations on noncompact
Riemannian manifolds}, Math. Meth. Appl. Sc. (2015)

\bibitem{[WeW09]}{G. F. Wei, W. Wylie},
\emph{Comparison geometry for the Bakry–Émery Ricci tensor}.
J. Differ. Geom. 83 (2009), 377--405.

\bibitem{[Wu10b]} {J-Y. Wu}, \emph{ Gradient estimates for a nonlinear diffusioon equation on complete manifold Journal}. J. Part. Diff. Eq. 23(1), 68--79

\bibitem{[Wu14]} {J-Y. Wu},
\emph{$L^p$-Liouville theorems on complete smooth metric measure spaces}. Bull. Sci. Math. 138 (2014), 510--539.

\bibitem{[Wu15]} {J-Y. Wu},
\emph{Elliptic gradient estimates for a weighted heat equation
and applications}. Math. Z. (2015) .

\bibitem{[Wu16]} {J-Y. Wu},
\emph{Elliptic gradient estimates for a nonlinear heat equation and applications}. Nonl. Anal. Theory, Mathods and Appl. 151 (2017), 1--17.


\bibitem{[YY08]}{Y. Y. Yang}, 
			{\it Gradient estimates for a nonlinear parabolic equation on Riemannian manifolds}. Proc Amer.
Math. Soc, 136 (2008), 4095--4102,
			 Acta. Math. Sin., 26(B)(2010), 1177--1182. 

\bibitem{[Y75]} {S-T. Yau}, 
\emph{Harmonic functions on complete Riemannian manifolds}. 
Commun. Pure Appl. Math. 28 (1975), 201--228.

			 
\bibitem{[Zhu11]} {X. Zhu}, 
\emph{Gradient estimates and Liouville theorems for nonlinear parabolic
equations on noncompact Riemannian manifolds}.
 Nonl. Anal., 74 (2011), 5141--5146.

\bibitem{[Zhu13]} {X. Zhu}, 
\emph{Hamilton’s gradient estimates and Liouville theorems for porous medium equations on noncompact Riemannian manifolds}. 
J. Math. Anal. Appl., 402(1) (2013), 201--206.


\bibitem{[Zhu16]} {X. Zhu}, 
\emph{Gradient estimates and Liouville theorems for linear and  nonlinear parabolic equations on Riemannian manifolds}. 
Acta Math Sc. 36(2)(2016), 514--526.
\end{thebibliography}
\end{document}